\documentclass{amsart}

\usepackage{amsmath,amscd,amssymb,amsthm,amsfonts}
\usepackage{enumerate}
\usepackage[]{graphicx}
\usepackage[all]{xy}
\usepackage{tensor}
\usepackage{tikz-cd}

\usepackage{xcolor}
\definecolor{red}{rgb}{0.8,0,0.1}
\definecolor{blue}{rgb}{0.1,0,0.8}

\newtheorem{theorem}{Theorem}[section]
\newtheorem{proposition}[theorem]{Proposition}

\newtheorem{lemma}[theorem]{Lemma}

\theoremstyle{definition}
\newtheorem{definition}[theorem]{Definition}
\newtheorem{example}[theorem]{Example}

\theoremstyle{remark}
\newtheorem{remark}[theorem]{Remark}

\numberwithin{equation}{section}

\def\C{{\mathcal C}}

\def\F{{\mathcal F}}
\def\G{{\mathcal G}}

\def\P{{\mathcal P}}

\def\to{\rightarrow}
\def\tto{\rightrightarrows}

\def\to{\rightarrow}
\def\tto{\rightrightarrows}

\newcommand{\hor}[0]{\mathrm{hor}}

\newcommand{\Curv}[0]{\mathrm{K}}

\newcommand{\curv}[0]{\mathrm{R}}
\newcommand{\mc}[0]{\mathfrak{mc}}

\newcommand{\pr}[0]{\mathrm{pr}}
\newcommand{\pp}[0]{\pi}
\newcommand{\qq}[0]{\mu}

\newcommand{\Hom}{\mathrm{Hom}}

\usepackage[foot]{amsaddr}

\begin{document}

\title{A note on Chern-Weil classes of Cartan connections}
\author{Luca Accornero$^1$}
\address{$^{1,3}$ Instituto de Matemática e Estatística - Universidade de São Paulo, Rua do Matão, 1010 - Butantã, São Paulo - SP, 05508-090, Brazil}
\email{$^{1}$luca.accornero@ime.usp.br, $^2$mateus.melo@ufes.br, 
$^3$ivanstru@ime.usp.br}
\author{Mateus M. de Melo$^2$}
\address{$^{2}$ Universidade Federal do Espírito Santo, Departemento de Matem\'atica, Avenida Fernando Ferrari 514 - Vitória, Espírito Santo - 29075-910, Brazil}
\author{Ivan Struchiner$^3$}
\thanks{L.A. was supported by the São Paulo Research Foundation (FAPESP), Brasil, Process Number 2024/22841-1. M.dM. acknowledges FAPESP for supporting his research through the grant 2019/14777-3. I.S. acknowledges FAPESP for supporting his research through the grant 2022/16310-8.}
\maketitle

\begin{abstract}
We present a construction of Chern–Weil characteristic classes for pairs of Cartan geometries sharing the same underlying group representation.

More precisely, given a “model” Cartan geometry, we define a subalgebra
of polynomials on its Atiyah algebroid such that any other Cartan geometry with the same underlying group representation comes with a characteristic map defined on such subalgebra and taking values in the deRham cohomology of the base manifold. The characteristic map recovers the classical Chern–Weil map of a Cartan connection when the "model" Cartan geometry
arises from a Klein geometry.
\end{abstract}

\setcounter{secnumdepth}{3}
\setcounter{tocdepth}{3}

\section{Introduction}
This note revolves around the construction of cohomological invariants of Cartan geometries. 

Cartan geometries are typically introduced as geometric structures on manifolds which are modeled infinitesimally on homogeneous manifolds, also known as Klein geometries -- see e.g. \cite{SharpeBook,CapSlovakParabolic}. Such a description presupposes the choice of an underlying Klein geometry $(G,H)$ consisting of a Lie group $H$ and a Lie subgroup $G$ of $H$ -- or, in slightly higher generality, an infinitesimal Klein geometry $(G,\mathfrak{h})$, also known as a model geometry. 

This follows the historical development: Cartan geometries were introduced as curved versions of Klein geometries. On the other hand, the modern formulation of Cartan geometries in terms of principal bundles and equivariant parallelisms on them does not require a model Lie algebra $\mathfrak{h}$ to be chosen. In fact, a Cartan geometry can be defined as a principal $G$-bundle $P$ together with a one form $\theta\in \Omega^1(P,V)$ -- called a Cartan connection -- taking value in a representation $V$ of $G$ that extends the adjoint representation. A compatible choice of Lie bracket on $V$ is needed to define the curvature of $\theta$, which is precisely the obstruction for $\theta$ to induce a morphism of Lie algebras from vector fields $\mathfrak{X}(P)$ with the usual Lie bracket to $V$ with the chosen Lie algebra structure.

While the notion of curvature is undoubtedly a central part of the theory of Cartan connections, it is desirable to develop a language which is to an extent model independent -- more precisely, dependent only on the underlying data $(G,V)$ -- as we try to argue below.

On the one hand, in many cases there are several model geometries which serve as model for the same geometric structure. For example, the usual metric on the round sphere can be modeled both on $(\mathrm{SO}_2, \mathfrak{so}_3)$, yielding a flat Cartan geometry, or on $(\mathrm{SO}_2, \mathfrak{euc}_2)$ -- where $\mathfrak{euc}_2$ is the Lie algebra of the group of isometries of the plane -- which yields a Cartan geometry with non-zero curvature. 

On the other hand, it often makes sense to consider as "models" instances of Cartan geometries that are not flat with respect to any natural Klein geometry/model geometry. A prominent example is contact geometry as a first order geometric structure (which we briefly go over in Section $4$).

This short note fits within these lines. The focus are Chern-Weil characteristic classes for Cartan geometries. Classically, those were addressed in \cite{AleMichCartan}. The result is a characteristic map constructed from the curvature of a Cartan connection and defined on the Lie algebra of invariant polynomials of the model -- we recall the construction in Theorem \ref{thm:chern_weil_CG}.

Here, given a "model" Cartan geometry $(Q,\omega)$ with underlying data $(G,V)$, we present a construction of Chern-Weil classes for any other Cartan geometry $(P,\theta)$ with the same underlying data. The resulting invariants should be understood as a measure how $(P,\theta)$ is different from $(Q,\omega)$. The main result is Theorem \ref{thm:main_thm} (followed by Example \ref{exm:grpd_CW_for_Klein_pairs}), which generalizes the work of \cite{AleMichCartan} while indicating a more intrinsic approach for obtaining invariants of geometric structures. 

\begin{theorem}
    Let $(Q,\omega)$ and $(P,\theta)$ be Cartan geometries over $M$ and $B$ respectively, with the same underlying data. Then, there exists a sub-algebra $S_\mathrm{char}^\bullet(Q,\omega)$ of the algebra 
    $S^\bullet\Gamma((TQ/G)^*)$ of polynomials on the vector bundle $TP/G$
    together with a characteristic map
    \[
    \kappa_P: S_\mathrm{char}^\bullet(Q,\omega) \to H^\bullet_\mathrm{dR}(B)
    \]
    which coincides with the classical characteristic map when $(Q,\omega)$ is the Cartan geometry associated to a Klein geometry.
\end{theorem}
The idea is the following: when having a Klein geometry $(G,H)$ at hand, one can interpret $H$ as the group of automorphisms of the Cartan geometry associated to the pair -- whose Cartan connection is simply the Maurer-Cartan form $\omega_H^\mathrm{MC}$ on $H$. Then $S_\mathrm{char}^\bullet(H,\omega_H^\mathrm{MC})$ is simply the algebra $S^\bullet_H(\mathfrak{h})$ of invariant polynomials on the Lie algebra of such group of automorphisms -- which, indeed, is a subalgebra of $\mathfrak{X}_G(H) \cong \Gamma(TH/G)$. 

For a general Cartan geometry $(Q,\omega)$, the group of automorphisms is not the appropriate notion -- for instance, it might contain only the identity. The way to proceed is to look at the gauge groupoid of the $Q$ instead -- a Lie groupoid over $M$ that encodes the principal bundle automorphisms of $Q$. The Lie algebra $\mathfrak{X}_G(Q)$ is nothing but the Lie algebra of sections of the associated Lie algebroid. The structure induced by $\omega$ on the gauge groupoid -- which has been addressed in \cite{CattafiCartanMult} -- is what allows one to introduce the appropriate generalization of $S^\bullet_H(\mathfrak{h})$.

We now offer some words on the methods. We make use of the theory of Lie groupoids and Cartan connections over them. Cartan connections are a very special instance of multiplicative one forms on Lie groupoids; those were studied in great detail in \cite{SalazarThesis,CrSalStr2012}, from which we borrow some results. Moreover, we refer more than once to \cite{CattafiCartanMult}, which investigates the gauge groupoid of a Cartan geometry. We also point out that the idea of using Lie groupoids and Lie algebroids to better understand as well as generalize Cartan geometries appears in the work of A.D. Blaom, e.g. \cite{BlaomInfSym}.

Finally, the reader familiar with groupoids and multiplicative structures will notice that the tools that we make use of are available in far greater generality than gauge groupoids and Cartan geometries. Indeed, 
our results can be appropriately stated and proved in such greater generality. The construction of $S_\mathrm{char}^\bullet(Q,\omega)$ came to us as a special case of a generalization of the Cartan model for equivariant cohomology to pairs consisting of a Lie groupoid together with a Cartan connection. Our results in this direction -- of which the present note can be seen as an application -- will appear in future work. 

The note is organized as follows. In Section $2$ we recall the basics on Cartan geometries and the construction of characteristic classes from \cite{AleMichCartan} -- Theorem \ref{thm:chern_weil_CG}. In Section $3$ we start by discussing gauge groupoids of Cartan geometries and recall the results of \cite{CattafiCartanMult}. Then, we explain how to present a Cartan geometry in terms of principal groupoid bundles with connections and, finally, we discuss our characteristic map and state our main result -- Theorem \ref{thm:main_thm} (together with Example \ref{exm:grpd_CW_for_Klein_pairs}). In Section $4$ we briefly discuss two prototypical example.  

We choose the Lie bracket on the Lie algebra of a Lie group and on the space of sections of the Lie algebroid of a Lie groupoid to be defined in terms of right invariant vector fields. All the actions and principal bundles are left actions unless otherwise specified.

\section{Cartan geometries and their characteristic classes}

\begin{definition}
    An \textbf{admissible data} for a Cartan geometry is a pair $(G,V)$ where $G$ is a Lie group and $V$ is a representation of $G$ that extends the adjoint representation on $\mathfrak{g}$ -- i.e. $\mathfrak{g}:=\mathrm{Lie}(G)$ is a sub-representation of $V$.
        
    A \textbf{model geometry} is an admissible data $(G,\mathfrak{h})$ such that $\mathfrak{h}$ is a Lie algebra and $\mathfrak{g}$ is a Lie subalgebra of $\mathfrak{h}$.
    
    A \textbf{Klein geometry} is a pair of Lie groups $(G,H)$ such that $G\hookrightarrow H$ is an embedded Lie subgroup -- thus, $(G,\mathfrak{h})$ is a model geometry when equipped with the restriction of the adjoint representation of $H$ to $G$. 
\end{definition}

\begin{definition}
    Let $(G,V)$ be an admissible data.
    A \textbf{Cartan connection} with underlying data $(G,V)$ on a principal $G$-bundle $P$ over $B$ is a one form $\theta\in \Omega^1(P,V)$ such that
    \begin{itemize}
        \item $\theta$ is a pointwise isomorphism -- hence, a parallelism on $P$;
        \item $\theta$ is compatible with the infinitesimal action, i.e. $\theta(a^P(\alpha))=\alpha$ for all $\alpha \in \mathfrak{g}$;
        \item $\theta$ is $G$-equivariant, i.e.
        \[
        L_g^*\theta = g\cdot \theta.
        \]
    \end{itemize}    
    The pair $(P,\theta)$ is called a \textbf{Cartan geometry}. A \textbf{morphism of Cartan geometries}
    \[\Phi:(P,\theta) \to (P',\theta')\]
    with underlying data $(G,V)$ 
    is a morphism of principal bundles such that
    \[
    \Phi^*(\theta') = \theta.
    \]
    
    Slightly more in general, if $(\phi,\psi):(G,V) \to (G',V')$ is an isomorphism of group representations, $(P,\theta)$ is a Cartan geometry with undelrlying data $(G,V)$, and $(P',\theta')$ is a Cartan geometry with undelrying data $(G',V')$, we say that a morphism of principal bundles $\Phi$ with underlying group morphism $\phi$ is a morphism between $(P,\theta)$ and $(P',\theta')$ with \textbf{underlying morphism of representations} $(\phi,\psi)$ when
    \[
    \Phi^*(\theta') = \psi\circ\theta.
    \]
    
    The notions of \textbf{automorphism} and \textbf{local morphism}/\textbf{local
    automorphism} are defined similarly.
\end{definition}
\begin{example}[Klein geometries]
    Klein geometries correspond to a specific class of Cartan geometries, presented by homogeneous spaces. The principal bundle is just $H$ itself together with the left action of $G$, and $B=H/G$ is the corresponding left coset space. The Cartan connection is simply the Maurer-Cartan form $\omega^{H}_\mathrm{MC}$ of $H$. Recall that the \textbf{Maurer-Cartan equation}
    \[
    d\omega^{H}_\mathrm{MC}+\frac{1}{2}[\omega^{H}_\mathrm{MC},\omega^{H}_\mathrm{MC}] = 0
    \]
    is satisfied.

   Since the Cartan connection is the Maurer-Cartan form of $H$, in this case the group of automorphisms can be identified with $H$. 
\end{example}
The definition of Cartan geometry is usually given in terms of a model geometry. One key reason is to be able to introduce the following notion.
\begin{definition}
    The \textbf{curvature} of a Cartan connection $\theta\in \Omega^1(P,\mathfrak{h})$ is the two form
    \[
    \mc^\theta_{\mathfrak{h}} = d\theta+\frac{1}{2}[\theta,\theta]\in \Omega^2(P,\mathfrak{h}). 
    \]
    A Cartan geometry with vanishing curvature is called \textbf{flat}.
\end{definition}
Similarly to what happens for principal connections, one has
\begin{lemma}\label{lem:basic_curvature}
    The two form $\mathfrak{mc}^\theta_\mathfrak{h}$ is basic -- i.e. there exists a vector bundle valued two form $\Curv^\theta_\mathfrak{h}\in \Omega^2(B,P[\mathfrak{h}])$, taking values in the associated bundle 
    \[
    P[\mathfrak{h}] := (P \times \mathfrak{h})/G
    \] 
    whose pullback to $P$ is $\mathfrak{mc}^\theta_\mathfrak{h}$. 
\end{lemma}


A classical result states that a flat Cartan geometry $(P,\theta)$ is locally isomorphic to a Klein geometry.
In particular, if $(P,\theta)$ is flat, then $P$ carries a local Lie group structure and $\theta$ can be interpreted as the corresponding Maurer-Cartan form. Thus, the theory of Cartan geometries can be interpreted as Lie theory with curvature \cite{SharpeBook}.

One more interpretation of Cartan connection is as restriction of principal ones. The following is proven e.g. in \cite{AleMichCartan}.
\begin{proposition}\label{prop:CG_are_principal_connections}
    Let $(G,H)$ be a Klein geometry and $(P,\theta)$ be a Cartan geometry on $B$ with underlying model $(G,\mathfrak{h})$. Then, the one form on $H\times P$
    \[
    \pr_1^*\omega_\mathrm{MC}^H - \pr_2^*\theta
    \]
    descends to a principal connection one-form $\Theta\in \Omega^1(P_H,\mathfrak{h})$ on the principal $H$-bundle
    \[
    P_H = (H\times P)/G\to B,
    \]
    where $G$ acts diagonally and the $H$-action on the quotient is induced by the left multiplication on the first factor. Moreover, the embedding
    \[
    P\to H\times P
    \]
    induces an inclusion of principal bundles with underlying group morphism given by the inclusion $G\subset H$, and the pull-back of $\Theta$ and its curvature $2$-form under such embedding are $\theta$ and $\mathfrak{mc}^\theta_\mathfrak{h}$ respectively.
\end{proposition}

Finally, the most relevant point for this note: the curvature of a Cartan connection gives rise to the characteristic classes of $P$. We refer to \cite{AleMichCartan} for a proof. 
\begin{theorem}\label{thm:chern_weil_CG}
    Let $(P,\theta)$ be a Cartan connection over $B$ with model geometry $(G,\mathfrak{h})$. Pairing with $\mc_\mathfrak{h}^\theta$ induces a characteristic map
    \[
    \kappa: S^\bullet_\mathfrak{h}(\mathfrak{h}^*)\to H^*_\mathrm{dR}(B)
    \]
    which is independent on the Cartan connection $\theta$ (here, $S^\bullet_\mathfrak{h}(\mathfrak{h}^*)$ denotes the algebra of invariant polynomials of $\mathfrak{h}$). In fact, the map is precisely the Chern-Weil map of $(P_H,\Theta)$ from Proposition \ref{prop:CG_are_principal_connections}. The map factors through the Chern-Weil map of the principal $G$-bundle
    \[
    \kappa: S^\bullet_\mathfrak{g}(\mathfrak{g}^*)\to H^*_\mathrm{dR}(B)
    \]
\end{theorem}

\subsection{Some model independent structure}

 Our first observation is that any Cartan geometry $(P,\theta)$ with underlying data $(G,V)$ comes with an induced skew symmetric $C^\infty(M)$-linear bracket on the space of sections of the bundle
 \[
 P[V] = P\times V/G,
 \]
where, as usual, the quotient with respect to the diagonal action. Sections of $P[V]$ can be identified with invariant functions $f\in C^\infty_G(P,V)$ and, applying the inverse of $\theta$, with sections of $TP/G$ -- i.e. invariant vector fields $X_f := \theta^{-1}(f) \in \mathfrak{X}_G(P)$. Keeping this in mind, the bracket is simply the differential of $\theta$:
 \begin{equation}\label{eq:fiberwise_bracket}
 \{\ ,\ \}: \Gamma(P[V]) \times \Gamma(P[V]) \to \Gamma(P[V]),\quad (f,g) \to d\theta(X_f,X_g),
 \end{equation}
where we observe that the equivariance of $\theta$ -- and hence of its differential -- ensures that the result is indeed an element of $\Gamma(P[V])$.
 

Via $\{\ ,\ \}$ one can extend the representation of $\mathfrak{g}$ on $V$ in the sense of the following Lemma.
\begin{lemma}\label{lem:ext_adjoint}
  For any function $f_\alpha\in C^\infty_G(P,\mathfrak{g})$ such that $f(p) = \alpha$ and any $f_v\in C^\infty_G(P,V)$ such that $f(p) = v$, it holds
  \[
  \{f_\alpha,f_v\}(p) = \alpha\cdot v,
  \]
  where $\cdot$ on the right hand side denotes the representation of $\mathfrak{g}$ on $V$. 
\end{lemma}
\begin{proof}
By deriving $g\cdot f_v = f_v\circ L_g$ at the identity in $G$ we get that
    \[
    \alpha\cdot v = \mathcal{L}_{a^P(\alpha)} f_v(p).
    \]
    On the other hand 
    \[
    \{f_\alpha,f_v\}(p) =d\theta(a^P\circ f_\alpha,X_{f_v})(p)
    \]
    where 
    $X_{f_v}\in \mathfrak{X}_G(P)$ is the invariant vector field $\theta^{-1}\circ f_v$. 
    
    Then, the claim follows from the identity
    \[
    \mathcal{L}_{X_{f_v}} f_\alpha = - \theta[a^P\circ f_\alpha,X_{f_v}],
    \]
    which is verified by using the definitions of the Lie derivative and the Lie bracket of vector fields in terms of flows.
\end{proof}
By exploiting $\{ \cdot , \cdot \}$, one also discovers a connection
\[
d^\theta: \Omega^\bullet(B, P[V]) \to \Omega^{\bullet+1}(B, P[V])
\]
 on the associated bundle $P[V]$.
 The derivation $d^\theta$ is defined on sections of $f\in \Gamma(P[V])$ by setting
\begin{equation}\label{eq:universal_connection}
d^\theta (f) (X)=  \mathcal{L}_{\hat{X}} f - \{\theta(\hat{X}),f\}
\end{equation}
where $\hat{X}$ is any $G$-invariant lift of $X\in \mathfrak{X}(B)$; the result does not depend on the chosen lift because of the identity
\[
\mathcal{L}_{a^P(\alpha)} f = \{\alpha, f\}
\]
which follows from Lemma \ref{lem:ext_adjoint} above and the equivariance of $f\in C^\infty_G(P,V)$. The curvature of $d^\theta$ will be denoted by $\curv^{d^\theta}\in \Gamma (\mathrm{Hom}(P[V],\Lambda^2 T^*B\otimes P[V]))$.

\section{Cartan geometries in terms of Lie groupoids}

\subsection{The gauge groupoid of a Cartan geometry}

Let $(P,\theta)$ be a flat Cartan geometry on $B$ with underlying Klein geometry $(G,H)$. The corresponding principal bundle with connection $(P_H,\theta^H)$ -- see Proposition \ref{prop:CG_are_principal_connections} -- is equipped with an $H$-equivariant map
\[\mu: P_H\to H/G,\quad [h,p]\to [h].\] Such a pair $(P_H,\mu)$ is equivalently encoded into a principal bundle $P_H \to B$ with structure groupoid the action groupoid  $H\ltimes H/G \rightrightarrows H/G$ and moment map $\mu: P_H \to H/G$. Moreover, the equivariance of $\theta^H$ with respect to $H$ can be understood as equivariance with respect to a certain action of the action groupoid on $TP_H$.

In this section, we make these remarks precise and extend them to Cartan geometries in a model independent way. The starting remark is that the action groupoid $H\ltimes H/G$ encodes the group of automorphism of the underlying Klein geometry. In fact, those coincide precisely with multiplication by elements of $H$.

In general, given a principal bundle, its automorphisms are encoded in a transtive groupoid over its base, called the \textbf{gauge groupoid}. Transitivity here is a reflection of the fact that given any two points $x,y\in M$ there exists a (local) automorphism of $P$ sending the fiber over $x$ to the fiber over $y$. The groupoid is defined by
\[
\G_P := (P\times P) / G \rightrightarrows P/G \cong M
\]
where the action is the diagonal action. An element of $\G_P$ is an equivalence class of pairs $(p,q)\in P\times P$ and is denoted by $[p:q]$. Thus, $[p:q]=[p':q']$ if and only if $gp=p'$ and $gq=q'$ for some $g\in G$; one deduces that bisections of $\G_P$ correspond to automorphisms of $P$. 

The following was proven in \cite{CattafiCartanMult}.
\begin{proposition}\label{prop:gauge_construction}
    Let $(P,\theta)$ be a Cartan geometry with underlying data $(G,V)$. Then 
    \begin{itemize}
        \item the associated vector bundle $P[V]$ is naturally a representation of $\mathcal{G}_P$, where for $[p:q]\in \G_P$ and $[(q,v)]\in P[V]$,
        \[
        [p:q] [(q,v)] = [(p,v)] \in P[V];
        \]
    \item the following one form
    \[
    (v_p,v_q) \mapsto \theta(v_p)-\theta(v_q)
    \]
    descends to the gauge groupoid $\mathcal{G}_P$ as a multiplicative form $\Theta\in \Omega^1(\mathcal{G}_P,t^*P[V])$ -- i.e. \begin{equation}\label{eq:multiplicativity}
    \Theta(dm(V_g,V_h)) = \Theta(V_g) + g\cdot \Theta(V_h)
    \end{equation}
    for any composable pair $(g,h)$, where $m$ is the groupoid multiplication and $g$ acts via the representation from the previous item;
    \item $\C_\Theta = \ker(\Theta)\subset T\mathcal{G}_P$ is complementary to both $T^s\G_P$ and $T^t\G_P$.
    \end{itemize}
    Moreover, local automorphisms of $(P,\theta)$ are in one to one correspondence with local bisections of $\G_P$ whose image in $\G_P$ is tangent to $\ker(\Theta)$.
\end{proposition}
The Lie algebroid of the gauge groupoid $\G_P$ is called the \textbf{Atiyah algebroid} of $P$ and can be identified with $TP/G\to M$ with anchor induced by the differential of the projection $P\to M$ and bracket on sections given by the bracket of $G$-invariant vector fields.
Thus, the isomorphism 
\[
P[V] \cong TP/G
\]
induced by $\theta$ can be read as an isomorphism between $P[V]$ and the Atiyah algebroid of $P$. Under this identification, the three items in the statement above can be packed together by saying that $\Theta\in \Omega(\G_P,t^*TP/G)$ is a \textbf{multiplicative Cartan connection} on the groupoid $\G_P$.
Notice that the representation on $TP/G$ resulting from the identification with $P[V]$ is simply
\[
w \mapsto \theta^{-1}_p\circ \theta_q(w),\quad ([p:q],w_{\mu(q)})\in \G_Q\times_M TP/G
\]
This interpretation makes also clear that there is a representation on $TM$
\[
(g,v)\in \G_Q\times_M TM \mapsto d\pi(g\cdot \hat{v})
\]
where $\hat{v}$ is any lift of $v$ to $TP/G$.
\begin{example}\label{ex:Klein_gauge_is_action}
    For a Klein geometry $(G,H)$,  $\G_H=H\times_G H \cong H\ltimes H/G $, where $H\ltimes H/G$ denotes the groupoid associated to the left action of $H$ on $H/G$ defined by
    \[
    (h,[k]) \to [kh^{-1}],\quad h\in H, [k]\in H/G;
    \]
    the isomorphism $\G_H=H\times_G H \cong H\ltimes H/G $ is
    \[
    [h:k]\in H\times_G H \mapsto (h^{-1}k,\pr(k))\in H\ltimes H/G.
    \]
    The multiplicative Cartan connection on $\G_H$ can be identified with the pullback to $H\ltimes H/G$ of the Maurer Cartan form on $H$; hence, bisections tangent to its kernel correspond to elements of $H$, as expected.

    At the level of Lie algebroids, the isomorphism $\G_H=H\times_G H \cong H\ltimes H/G $ described above induces the following isomorphism
    \[
    [(h,\alpha)]\in H[h] \mapsto (h^{-1}\cdot \alpha,\pr(h)) \in h\ltimes H/G.
    \]
\end{example}

Multiplicative Cartan connections interact with the Lie groupoid structure nicely. One instance of such interaction is that a multiplicative connection on a Lie groupoid $\G$ induces, by differentiation, a connection on the corresponding algebroid $A$ that satisfies an infinitesimal version of multiplicativity. We refer to \cite{SalazarThesis,CrSalStr2012} for details and limit ourselves to give the formula which is relevant for our purposes: the connection
\[
D: \Omega^\bullet(M,TP/G) \to \Omega^{\bullet+1}(M,TP/G)
\]
resulting from differentiation is defined on sections $V\in \Gamma(TP/G) \cong \mathfrak{X}_G(P) $ by
\[
D_X V = \nabla_{V} \hat{X} + [\hat{X},V]
\]
where $X\in \mathfrak{X}(M)$, $\hat{X}\in \mathfrak{X}_G(P)$ is a lift of $X$, and $\nabla$ denotes the representation of $TP/G$ on itself obtained by differentiating the representation in the first item of Proposition \ref{prop:gauge_construction}. 

The equivariant function $f_{\nabla_{V} \hat{X}}$ corresponding to $\nabla_{V} \hat{X}$ is given by
\[
\mathcal{L}_{V} \theta(\hat{X}).
\]
Hence, under the identification $\mathfrak{X}_G(P)\cong \Gamma(P[V])$ provided by $\theta$, $D$ corresponds precisely to the connection
\[
d^\theta: \Omega^\bullet(M,P[V]) \to \Omega^{\bullet+1}(M,P[V])
\]
induced by $\theta$ and discussed in the previous section. 

\subsection{Gauge groupoids as generalized model geometries}

Let $(Q,\omega)$ be a Cartan geometry over $M$ with underlying data $(G,V)$. Throughout this section $(\G_Q,\Omega)$ denotes its gauge groupoid with the multiplicative Cartan connection constructed above. Recall that a \textbf{principal $\G_Q$-bundle} with momentum map $\mu:P\to M$ and projection $\pi: P\to B$, denoted by $Q\overset{\mu}{\leftarrow} \mathcal{P}\to B$, is an action of the Lie groupoid $\G_Q$ along the map $\mu$ such that the map $\pi$ is $\G_Q$-invariant and the action map
\[
m_P: \G_Q\times_M P \to P\times_B P
\]
is a diffeomorphism. Recall also that the infinitesimal action underlying a principal bundle is the vector bundle map
\[
a^P:\qq^*A \to TP,\quad \alpha \to dm(\alpha_{\qq(p)},0_p);
\]
this is anchored -- i.e. $\rho = d\mu\circ a^P$ -- and compatible with the Lie bracket on sections of $A$ -- i.e. $a^P([\qq^*\alpha,\qq^*\beta])=[a^P(\qq^*\alpha),a^P(\qq^*\beta)]$ for any pair of sections $\alpha,\beta\in \Gamma(A)$.

In general, an action of a Lie groupoid on $P$ does not induce an action on $TP$; as a consequence, the notion of principal connection for principal groupoid bundles is non-trivial. However, the multiplicative Cartan connection $\Omega$ on $\G_Q$ allows to lift the action to $TP$ as follows
\[
(g,v_p) \in \G_Q\times_M TP \to dm_P(\hor^\Omega_g(d\mu(v_p)),v_p)
\]
where $\hor^\Omega_g(d\mu(v_p))$ is the only tangent vector in $T_g\G_Q$ that projects to $d\mu(v_p)$ via $ds$ and is tangent to $\ker(\Omega)$. The fact that the formula above defines an action is a consequence of the multiplcativity of $\Omega$.

\begin{definition}\label{def:princ_conn}
    A \textbf{principal $(\G_Q,\Omega)$-connection} one form on a principal $\G_Q$-bundle $Q\overset{\mu}{\leftarrow} \mathcal{P}\to B$ consists of a one form $\Theta\in \Omega^1(\P,\mu^*A)$ which is pointwise surjective and satisfying the multiplicativity equation
    \[
    \Theta(dm(V_g,V_p)) = \Theta(V_g) + g\cdot \Theta(V_p)
    \]
    for all composable vectors $V_g\in T_gG$ and $V_p\in T_pP$. The pair $(\P,\Theta)$ is also called a \textbf{principal $(\G_Q,\Omega)$-bundle}. An \textbf{isomorphsm} between principal $(\G_Q,\Omega)$-bundles $(\P,\Theta)$ and $(\P',\Theta')$ is an isomorphism of principal $\G_Q$-bundles that pulls back $\Theta'$ to $\Theta$.
\end{definition}
\begin{example}\label{ex:std_PB}
The action of $(\G_Q,\Omega)$ on itself by left multiplication is an example of principal $(\G_Q,\omega)$-bundle -- called the \textbf{standard} principal $(\G_Q,\omega)$-bundle.
\end{example}
Notice that the kernel of a principal $(\G_Q,\Omega)$-connection one form $\Theta$ defines an Ehresmann connection with respect to $\pi:\P\to B$; moreover the differential of the action map induces a map
\[
\ker(\Omega)\times_{TM}\ker(\Theta) \to \ker(\Theta)\subset T\P_Q.
\]
We defer a more complete theory of principal connections -- which can be introduced in the context of groupoids equipped with multiplicative Cartan connections and do not rely on the specific structure of $\G_Q$ -- to later work. Here, we limit ourselves to the statements that are essential for our purposes, starting with the one below.
\begin{proposition}\label{thm:CG_are_gauge_PB}
    Let $(Q,\omega)\to M$ be a Cartan geometry with underlying data $(G,V)$. There is a one to one correspondence between isomorphism classes of Cartan geometries $(P,\theta)\to B$ with underlying data $(G,V)$ and isomorphism classes of principal $(\mathcal{G}_Q,\Omega)$-bundles $M\overset{\mu}{\leftarrow} (\mathcal{P}_Q,\Theta)\to B$.
\end{proposition}
\begin{proof}
    
    Starting from a Cartan geometry $(P,\theta)$ over $B$, we set
    \[
    \P_Q:= (Q\times P)/G,
    \]
    where, as usual, the quotient is with respect to the diagonal action. The structure of principal bundle over $B$ is induced by the canonical structure of right principal $Q\times Q$-bundle over $P$ on $Q\times P$, with momentum map given by the first projection and action given by
    \[
    (p,q_1)\cdot (q_1,q_2) = (p,q_2),\quad p\in P,\ q_1,q_2\in Q.
    \]
    Since all the relevant structure maps are $G$-invariant, it follows that $\G_Q$ acts on $\P_Q$ principally with base $B$. Notice that, to obtain a left principal bundle, we need to compose the action map with the inversion of $\G_Q$ -- which simply swaps the components of a representative. Concerning the form $\Theta\in \Omega^1(\P_Q,\mu^*A)$, it is obtained from the form (recall that $A\cong P[V])$
    \[
    \pr_Q^*\omega - \pr_P^*\theta\in \Omega^1(Q\times P, V),
    \]
    which vanishes along the $G$ orbits and is $G$-equivariant, and thus descends to the quotient under the $G$-action. The principality of the induced form $\Theta$ (Definition \ref{def:princ_conn}) is then readily checked by working at the level of the principal $Q\times Q$-bundle $Q\times P$. 

    If we start from a principal $(\G_Q,\Omega)$-bundle $(\P,\Theta)$, we construct a Cartan geometry $(P,\theta)$ by choosing a point $x\in M$ and restricting the action of $\G_Q$ to the action of its isotropy group $(\G_Q)_x$ at $x$ on $P_x=\mu^{-1}(x)$ equipped with the restriction $\theta_x:= \Theta|_{\mu^{-1}(x)}$ of $\Theta$. The multiplicativity equation \eqref{eq:multiplicativity} of $\Theta$ ensures that $(P_x,\theta_x)$ is indeed a Cartan geometry with underlying data $((\G_P)_x,A_{x})$. Different choices of base-point $x$ yield isomorphic geometries. Indeed, for $x,y\in M$, any arrow $[p:q]\G_P$ with source $x$ and target $y$ induces, by coniugation, isomorphism of Lie groups
    \[
    \mathcal{C}_{[p:q]}: (\G_P)_x \to (\G_P)_y, \quad g\in (\G_P)_x\mapsto \left( [p:q] \cdot g \cdot [q:p]\right) \in (\G_P)_y
    \]
    and a corresponding isomorphism of principal bundles from $\qq^{-1}(x)$ to $\qq^{-1}(y)$. The fact that the latter is compatible with the forms $\theta_x$ and $\theta_y$ follows once again from the multiplicativity equation \eqref{eq:multiplicativity}.

    Finally, direct checks show that isomorphic Cartan geometries yield isomorphic principal groupoid bundles with connection and viceversa, and the two construction are (up to isomorphisms) inverse of each other.
\end{proof}

\begin{example}\label{ex:action_principal_bundle_conn}
    Lets assume that $(Q,\omega)$ is the Klein geometry corresponding to the pair $(G,H)$. We already observed -- example \ref{ex:Klein_gauge_is_action} -- that 
    \[
    (H\ltimes H/G,\pr_1^*\omega^\mathrm{MC}_H) \cong (\G_Q,\omega)
    \]
    Here, we further point out that a principal $H\ltimes H/G$-bundle $\P$ is the same thing as a principal $H$-bundle with an $H$-equivariant map to $H/G$. For the bundle constructed above, the equivariant map is
    \[
    \P_H \to H/G, \quad (h,p) \to \pr(h). 
    \]
    Moreover, the connection one form $\Theta\in \Omega^1(\P_H,\mu^*(\mathfrak{h}\ltimes H/G))$ is in fact an Ehresmann connection for the principal $H$-bundle $\P_H\to B$. Thus, our Proposition \ref{thm:CG_are_gauge_PB} recovers the well known correspondence between Cartan geometries and principal connections that we recalled in Proposition \ref{prop:CG_are_principal_connections}.
\end{example}

\begin{example}
    In the particular case when $(P,\theta)=(Q,\omega)$, the construction above yields $(\G_Q,\omega)$ acting upon itself -- Example \ref{ex:std_PB}.
\end{example}

\subsection{Curvature and Chern-Weil construction}

Throughout the rest of this note, $(Q,\omega)$ and $(P,\theta)$ are Cartan geometries over $M$ and $B$ respectively, with the same underlying data $(G,V)$, and $M \overset{\mu}{\leftarrow}(\P_Q,\Theta)\overset{\pi}{\to} B$ is the principal $(\G_Q,\Omega)$-bundle corresponding to $(P,\theta)$ via Proposition \ref{thm:CG_are_gauge_PB}. 

\begin{definition}
    The \textbf{Maurer-Cartan} expression of $(\P_Q,\Theta)$ is the two form $\mc^\Theta \in \Omega^2(\P_Q,\mu^*Q[V])$ 
    defined as 
    \[\mc^\Theta = (d_\mu^\omega) \Theta+\frac{1}{2}\{\Theta,\Theta\}^\omega,\]
    where $d_\mu^\omega$ is the pullback via $\mu$ of the connection $d^\omega$ while $\{\ ,\ \}^\omega$ is the fiberwise bracket on $P[V]$ induced by $\omega$, see \eqref{eq:fiberwise_bracket}.
\end{definition} 
The first item in the proposition below is discussed in \cite{AccCr2023}, while the second is a direct check.
\begin{proposition}\label{prop:properties_of_mc}
    The form $\mc^\Theta$ is 
    \begin{itemize}
        \item horizontal -- i.e. $\mc^\Theta(a^P(\alpha_{\qq(p)}),\cdot) = 0$ for all $\alpha_{\qq(p)}\in \qq^*A$;
        \item multiplicative with respect to the Maurer-Cartan expression $\mc^\Omega$ of the standard principal $(\G_Q,\Omega)$-bundle from Example \ref{ex:std_PB}, i.e. for all $(g,p)\in G\times_{M} P$
        \[
        (\mc^\Theta)_{gp}\circ (dm_{\P_Q},dm_{\P_Q}) = (\mc^\Omega)_{g} \circ (d\pr_1,d\pr_1) + g\cdot(\mc^\Theta)_{p}\circ(d\pr_2,d\pr_2). 
        \]
    \end{itemize}
\end{proposition}

Notice that, from the first item in the proposition above, it follows
    \[
    \mc^\Theta(X,Y) = d^\omega_\mu \Theta(X^\Theta,Y^\Theta) = -\Theta([X^\Theta,Y^\Theta]), \quad X,Y \in \mathfrak{X}(\P_Q)
    \]
    where $X^\Theta$ and $Y^\Theta$ denote the $\ker(\Theta)$-components of $X$ and $Y$ with respect to the splitting $T\P_Q = \ker(\Theta)\oplus T^\pi(\P_Q)$.

In turn, this shows that $\mc^\Theta$ can be interpreted as a section
\[
\mc^\Theta : P \to \mathrm{Hom}(\pp^*\Lambda^2 TB,\qq^*Q[V])
\]
which vanishes if and only if $\ker(\Theta)$ is involutive.

Moreover, by focusing on the case of the standard principal bundle $(\P_Q,\Theta)=(\G_Q,\Omega)$ with the action by left multiplication, we deduce that $\mc^\Omega$ is a groupoid cocycle with values in the representation $\mathrm{Hom}(\Lambda^2TM,Q[V])$ induced by the representation of $\G_Q$ on $P[V]$ and $TM$. This prompts one to investigate 
its linearization.
\begin{lemma}\label{lem:lineariz_curv}
    The Lie algebroid cocycle 
    \[\mathrm{Lie}(\mc^\Omega)\in \Gamma(\Hom(Q[V],\Hom(\Lambda^2TM,Q[V]))\]
    obtained from differentiation of $\mc^\Omega$ and using the identification $\mathrm{Lie}(\G_Q) = A \cong Q[V]$ is precisely the curvature of the connection $d^\theta$ from \eqref{eq:universal_connection}.
\end{lemma}
\begin{proof}
    As we already recalled $d^\theta$ can be interpreted as the  connection on the Lie algebroid $A$ obtained from $\Theta$ by linearization. The fact that linearization is compatible with taking curvature is then a general fact about multiplicative forms, discussed in \cite{SalazarThesis}.
\end{proof}

\begin{definition}
    We say that $(\P_Q,\Theta)$ is \textbf{flat} when $\mc^\Theta$ = 0. Similarly, we say that $(\G_Q,\Omega)$ is a \textbf{flat groupoid} if $\mc^\Omega = 0$. 
\end{definition}
\begin{example}\label{ex:Klein_pair_flat}
    Given a Klein geometry $(G,H)$ and setting $(Q,\omega) = (H,\omega^\mathrm{MC}_{H})$, from Example \ref{ex:Klein_gauge_is_action} we see that $\mc^\Omega = 0$.
\end{example}
Observe that a consequence of the second item in Proposition \ref{prop:properties_of_mc} is that if $(\P_Q,\Theta)$ is flat, then $(\G_Q,\Omega)$ is a flat groupoid. Characteristic classes for flat principal bundles $(\P,\Theta)$ were studied in \cite{AccCr2023}. On the other hand, we also deduce:
\begin{proposition}
    If $(\G_Q,\Omega)$ is a flat groupoid, then $\mc^\Theta$ is basic -- i.e. there exists a two form $\Curv^\Theta \in \Omega^2(B,P[Q[V]])$ such that 
    \[
    \mc^\Theta = \pi^*\Curv^\Theta.
    \]
\end{proposition}
If $(P,\theta)$ has an underlying Klein geometry $(G,H)$ and $(Q,\omega)=(H,\omega^\mathrm{MC}_{H})$ is the geometry corresponding to the pair,
$\Curv^\Theta$ is nothing but the form corresponding to the curvature of the principal connection from Proposition \ref{prop:CG_are_principal_connections} -- recall also Example \ref{ex:action_principal_bundle_conn}.




The bundle $(\P_Q,\Theta)$ and the form $\mc^\Theta$ introduced above provide a way to introduce Chern-Weil like invariants of Cartan geometries that make reference to the Cartan geometry $(Q,\omega)$ rather than to a Klein geometry. 

Let $S^k(TQ/G)^*$ denote $k$-th symmetric product of $TQ/G$ with itself, and by $S^\bullet((TP/G)^*)$ the corresponding algebra. The horizontality of $\mc^\Theta$ ensures that the map
\[
\kappa_\Theta: \Gamma(S^\bullet(TQ/G)^*)) \to \Omega^{2\cdot \bullet}(\P_Q),\quad \gamma\to \gamma(\mc^\Theta,\ldots, \mc^\Theta)
\]
takes values in the sub-algebra of horizontal forms. Here, and throughout this subsection, we identify $\gamma\in \Gamma(S^\bullet(TQ/G)^*))$ with the corresponding section of the pullback bundle $\qq^*S^\bullet(TQ/G)^*$.

We first investigate which elements of $\Gamma(S^\bullet(TQ/G)^*)$ take value in the smaller sub-algebra basic forms. Recall that a horizontal form $\tau$ on the principal $\G_Q$-bundle $\P_Q$ is basic if and only if
\[
\tau\circ (dm_{\P_Q},\dots ,dm_{\P_Q})|_{(g,p)} =  \tau\circ \pr_{\P_Q}|_{(g,p)},\quad (g,p)\in G_Q\times_M \P_Q.
\]

\begin{lemma}
    Let $\Gamma(S^\bullet(TQ/G)^*)_\mathrm{bas}$ be the subspace of sections of $\Gamma(S^\bullet(TQ/G)^*)$ that are
    \begin{itemize}
    \item invariant under the representation of $\G_Q$ on $S^\bullet(TQ/G)^*$ -- i.e. for all $g\in \G_P$ it holds
    \[
    g\cdot \gamma_{s(g)} = \gamma_{t(g)};
    \]
    \item annihilated by $\mc^\Omega$ -- i.e. if $\gamma \in \Gamma(S^\bullet(TQ/G)^*)_\mathrm{bas}$ then
    \[
    \gamma\circ (\mc^\Omega,-,\ldots,
    -) = 0 \in \Omega^2(\G_Q)\otimes \Gamma(S^{\bullet -1}((TQ/G)^*))
    \]
\end{itemize}
The restriction of $\kappa_\Theta$ to $\Gamma(S^\bullet(TQ/G)^*)_\mathrm{bas}$ takes values into basic forms on $P$.
\end{lemma}
\begin{proof}
    We need to show that, if $\gamma\in \Gamma(S^\bullet(TQ/G)^*)_\mathrm{bas}$, then the form
    \[
    \kappa_\Theta(\gamma) := \gamma \circ (\mc^\Theta,\ldots,\mc^\Theta)\in \Omega^{2\bullet}(P)
    \]
    is basic. 

    The second item in proposition \ref{prop:properties_of_mc} yields, for all $(g,p)\in \G_Q\times_M\P_Q$

    \begin{align*}
    \kappa_\Theta(\gamma)\circ (dm_{\P_Q},\circ \ldots,dm_{\P_Q})|_{(g,p)} = \gamma\circ (\mc^\Omega\circ d\pr_{\G_Q} + g\cdot\mc^\Theta\circ d\pr_{\P_Q},\\
    \ldots, \mc^\Omega\circ d\pr_{\G_Q} + g\cdot\mc^\Theta\circ d\pr_{\P_Q})|_{(g,p)},
    \end{align*}
    which, since $\gamma$ is annihilated by $\mc^\Omega$, is equal to 
    \[
    \gamma \circ (g\cdot\mc^\Theta\circ \pr_{\P_Q},\ldots, g\cdot\mc^\Theta\circ \pr_{\P_Q})|_{(g,p)}.
    \]
    Then, the statement follows because $\gamma$ is invariant under the action of $\G_P$.
\end{proof}


Thus, we find out that $(\P_Q,\Theta)$ comes with a map, still denoted by $\kappa_\Theta$,
\[
\kappa_\Theta: \Gamma(S^\bullet(TQ/G)^*)_\mathrm{bas}\to \Omega^{2\bullet}(B).
\]
Now, let $\Gamma(S^\bullet(TQ/G)^*)_\mathrm{char}\subset \Gamma(S^\bullet(TQ/G)^*)_\mathrm{bas}$ be the subset consisting of polynomials that are also closed under the connection on $S^\bullet(TQ/G)^*$ induced by $d^\omega$. A direct computation shows the following:
\begin{proposition}
    The restriction of $\kappa_\Theta$ to $\Gamma(S^\bullet(TQ/G)^*)_\mathrm{char}$ takes values in the sub-algebra of closed forms.
\end{proposition}
\begin{proof}
     It is sufficient to prove that 
     \[
     \kappa_\Theta \circ d^\omega = d\circ \kappa_\Theta
     \]
     Due to the fact that $d$ and $d^\omega$ are derivations, it is sufficient to work on sections of $T^*Q/G$. Then
     \[
     \gamma(\mc^\Theta) = \langle \gamma,\mc^\Theta\rangle,\quad \gamma\in \Gamma(T^*P/G)
     \]
     and we have
     \[
     d(\langle \gamma,\mc^\Theta\rangle) =\langle d^\omega\gamma, \mc^\Theta \rangle + \langle \gamma, d^\omega_{\qq}\mc^\Theta \rangle.
     \]
     The first term vanishes due to our assumption that $\gamma$ is closed under $d^\omega$. Thus, we focus on the second term.
     
    Since $d(\langle \gamma,\mc^\Theta\rangle)$ is basic, it is sufficient to show that 
    \[
    \langle \gamma, d^\omega_{\qq}\mc^\Theta \rangle
    \]
    vanishes when evaluated along vectors tangent to $\ker(\Theta)$ . Moreover, as explained right after Proposition \ref{prop:properties_of_mc}, 
    \[
    \mc^\Theta(X,Y) = d^\omega_\mu \Theta(X^\Theta,Y^\Theta),\quad X,Y\in \mathfrak{X}(P)
    \]
    where $X^\Theta$ and $Y^\Theta$ denote the components of $X$ and $Y$ tangent to $\ker(\Theta)$. Thus, we are left with checking that
    \[
    (d^\omega_\mu)^2\Theta (X,Y,Z) = 0,\quad X,Y,Z\in \ker(\Theta),
    \]
    which is done by means of a direct calculation.
 \end{proof}

Thus, $\kappa_\Theta$ induces a map valued in the de Rham cohomology of $B$, for which we keep the same notation:
\[
\kappa_\Theta: \Gamma(S^\bullet(TQ/G)^*)_\mathrm{char} \to H^{2\cdot \bullet}_\mathrm{dR}(B). 
\]
As in classical Chern-Weil theory, the homotopy invariance of de Rham cohomology 
implies that the cohomology class of $\kappa_\Theta(\gamma)$, for $\gamma\in S^\bullet(TQ/G)^*$, does not depend on $\Theta$ -- i.e. any other choice of principal $(\G_Q,\Omega)$-connection induces the same map in cohomology. Hence, we can finally state
\begin{theorem}\label{thm:main_thm}
    Let $(Q,\omega)$ be a Cartan geometry with underlying data $(G,V)$. For any other Cartan geometry $(P,\theta)$ with the same underlying data, composition with the Maurer-Cartan expression $\mc^\Theta$ induces a characteristic map
    \[
    \kappa_P: S^\bullet_\mathrm{char}(Q[V])\to H^\bullet_\mathrm{dR}(M),
    \]
    which does not depend on the choice of $\theta$ among Cartan connections on $P$ with underlying data $(G,V)$.
\end{theorem}


\begin{example}\label{exm:grpd_CW_for_Klein_pairs}
    Let us go back to Klein geometries $(G,H)$ -- i.e. let us assume that $(P,\theta)$ is modelled on $(G,H)$ and let's take $(Q,\omega)$ to be $H$ with its Maurer-Cartan form. Since $\G_H$ has vanishing Maurer-Cartan expression -- see Example \ref{ex:Klein_pair_flat} -- the structure of $S^\bullet_\mathrm{char}$ greatly simplifies. The vanishing of $\mc^\Omega$ implies the flatness of the connection $d^\omega$ -- whose flat sections correspond to elements of $\mathfrak{h}$. Moreover, from $(\G_H,\Omega)\cong (H\ltimes H/G,\pr_1^*\omega^\mathrm{MC}_H)$ -- Example \ref{ex:Klein_gauge_is_action} -- we see that the invariance under the groupoid action translates into the invariance under the adjoint action of $H$ on $\mathfrak{h}$. In other words
    \[
    \Gamma(S^\bullet(Q[V]))_\mathrm{char}\cong S^\bullet_H(\mathfrak{h}).
    \]

    Moreover, under the identification between principal $(H\ltimes H/G,\pr_1^*\omega_H^\mathrm{MC})$-bundles and principal $H$-bundles with connection recalled in example \ref{ex:action_principal_bundle_conn}, the map from Theorem \ref{thm:main_thm} becomes the classical Chern-Weil map -- thus recovering the construction from \cite{AleMichCartan}, see Theorem \ref{thm:chern_weil_CG}.

\end{example}

\begin{remark}
    While Theorem \ref{thm:chern_weil_CG} is stated in terms of model geometries rather than Klein geometries, we stress here that our construction can be formulated in similar infinitesimal terms -- so that it recovers Theorem \ref{thm:chern_weil_CG} when applied to Cartan geometries with a choice of underlying model. Indeed, in the definition of $\Gamma(S^\bullet(Q[V]))_\mathrm{char}$, one simply replaces the invariance under the representation of $\G_Q$ with invariance under the corresponding representation of its Lie algebroid $Q[V]$.
\end{remark}

\begin{remark}\label{rmk:classical_construction_is_Cartan_model}
    Going back to Example \ref{exm:grpd_CW_for_Klein_pairs}, a closer look reveals an insightful interaction with the Cartan model for equivariant cohomology that we will address, in greater generality, in upcoming work. Let us briefly mention what happens in the setting of this note. 
    
     In example \ref{ex:action_principal_bundle_conn} we observed that principal $(H\ltimes H/G,\pr_1^*\omega_H^\mathrm{MC})$-bundles correspond to principal $H$-bundles with connection together with an $H$-equivariant map $\qq$  to $H/G$. Such a map induces a map at the level of equivariant cohomology
    \[
    \qq: H^\bullet_H(H/G) \to H_H^\bullet(P) 
    \]
    For a compact group $H$, the equivariant de Rham theorem allows to model both cohomologies in terms of $S^\bullet_H(\mathfrak{h})$-valued invariant forms on $H/G$ and $P$ respectively -- the so called Cartan model for equivariant cohomology. For a non compact $H$, $\qq$ still induces a map at the level of Cartan models. 
    
    Furthermore $H^\bullet_H(H/G)\cong S^\bullet_H(\mathfrak{h})$ and $H_H^\bullet(P) \cong H^\bullet(B)$. The former isomorphism follows by using the Cartan model and the fact that the action is transitive. The latter isomorphism is due to principality and, at the level of the Cartan model, it can be constructed explicitly by using a principal coonnection and classical Chern-Weil theory. 
    
    To sum up, the cohomology map induced by $\qq$ can be interpreted as a map
    \[
    S^\bullet_H(\mathfrak{h}) \to H_\mathrm{dR}(B)
    \]
    and, as such, it coincides with the characteristic map of the Cartan geometry $(P,\theta)$ modeled on $(G,H)$. Our construction generalizes the map to the gauge groupoid with Cartan connection $(\G_Q,\Omega)$ of a possibly non-flat Cartan geometry $(Q,\omega)$; this can in fact be interpreted as a generalization of the Cartan model from action groupoids to pairs of the form $(\G_Q,\Omega)$. In upcoming work, we will show that such a generalization is possible in the context of groupoids equipped with a Cartan connection. 
\end{remark}

\section{Examples}
We now apply the results above to two prototypical examples: $G$-structures with connections, where a canonical model geometry (in fact a Klein geometry) is always available, and contact geometry, where no such model geometry exists.

\subsection*{G-structures with connections}
One class of examples is obtained by looking at $G$-structures with connections. This class admits a canonical underlying Klein geometry, and fits the classical story

Recall that a $G$-structure on a manifold $B$, for $G\hookrightarrow \mathrm{Gl}_n(\mathbb{R})$, is a reduction of the frame bundle of $B$
\[(P,\tau)\hookrightarrow (\mathrm{Fr}(B), \tau^\mathrm{taut}),\]
where $\tau \in \Omega^1(P,\mathbb{R}^n)$ and $\tau^\mathrm{taut}$ is the tautological or solder form
\[
\tau^{\mathrm{taut}}_\phi(v) \to \phi^{-1}(d\pi(v))\in \mathbb{R}^n\].

A $G$-structure $(P,\tau)$ can be turned into a Cartan geometry by choosing a principal connection $\omega\in \Omega^1(P,\mathfrak{g})$ and setting
\[
\theta:=\tau+\omega.
\]

\begin{remark}
In some cases there is a canonical choice of a connection making the $G$-structure into a Cartan geometry. This is the case, for example, of $\mathrm{O}_n(\mathbb{R})$- structures (Riemannian structures) or more generally holonomy reductions of Riemannian structures.
\end{remark}

One instance of $G$-structure is the \textbf{flat $G$-structure} on $\mathbb{R}^n$
\[
(P^\textrm{$G$-flat}, \tau):= (\mathbb{R}^n \times G, A\cdot \pr_{\mathbb{R}^n}^*\omega^\mathrm{MC}_{\mathbb{R}^n}) \hookrightarrow (\mathrm{Fr}(\mathbb{R}^n), \tau^\mathrm{taut})
\]
where $\omega^\mathrm{MC}_{\mathbb{R}^n}$ is the Maurer-Cartan form of $\mathbb{R}^n$ with its abelian Lie group structure and we use the notation $A\cdot \pr_{\mathbb{R}^n}^*\omega^\mathrm{MC}_{\mathbb{R}^n})$ for the form 
\[
v\in T_{x,A}(P^\textrm{$G$-flat}) \mapsto A\cdot \omega^\mathrm{MC}_{\mathbb{R}^n}(v).
\]
The $G$-structure $(\mathbb{R}^n \times G, \pr_1^*\omega^\mathrm{MC}_{\mathbb{R}^n})$ carries a canonical flat principal connection, and thus can be turned into a flat Cartan geometry. Indeed, the one-form
\[
\theta^\textrm{$G$-flat}:=\tau +\pr_{G}^*\omega^\mathrm{MC}_G
\]
is a flat Cartan connection on $P^\textrm{$G$-flat}$ modeled on the Klein geometry $(G, \mathrm{Aff}^n_G:= G\ltimes \mathbb{R}^n)$.

Performing the gauge construction on $(P^\textrm{$G$-flat},\theta^\textrm{$G$-flat})$ we discover the groupoid
\[
\G^\mathrm{flat}:=\mathbb{R}^n\times G\times \mathbb{R}^n \tto \mathbb{R}^n
\]
with source map $s(y, A,x) = x$, target map $s(y, A,x) = y$ and multiplication $(z,A,y)\cdot (y,B,x) = (z,AB,x)$ for all $(z,A,y)$ and $(y,B,x)$ in $\mathbb{R}^n\times G\times \mathbb{R}^n$. The resulting Cartan connection is
\[\Theta^\textrm{$G$-flat}:=\pr_{G}^*\omega^\mathrm{MC}_G+\delta\omega^\mathrm{MC}_{\mathbb{R}^n},\]
where $\delta\omega^\mathrm{MC}_{\mathbb{R}^n}|_{(y,A,x)} = \omega^\mathrm{MC}_{\mathbb{R}^n}|_{y} - A\cdot \omega^\mathrm{MC}_{\mathbb{R}^n}|_{x}$. As expected, there is an isomorphism of groupoids with Cartan connections
\[
(\G^\mathrm{flat},\Theta^\textrm{$G$-flat})\cong  (\mathrm{Aff}^n_G\ltimes \mathbb{R^n}, \pr_1^*\omega^\mathrm{MC}_{\mathrm{Aff}^n_G})
\]
defined by 
\[
(y,A,x) \mapsto ((A,y-Ax),x).
\]

Now, Theorem \ref{thm:CG_are_gauge_PB} guarantees that any $G$-structure with connection $(P,\theta)$ can be presented as a principal $(\G^\mathrm{flat},\Theta^\mathrm{flat})$-bundle $(\P,\Theta)$ -- i.e. a principal $\mathrm{Aff}_G$-bundle $\P$ with a principal connection $\Theta$. 

This construction is within the scope of the classical Theorem \ref{thm:chern_weil_CG}. Thus, the Chern-Weil classes of $(\P,\Theta)$ correspond to those of $(P,\theta)$ and the characteristic map is defined on $S^\bullet_{\mathfrak{aff}_\mathfrak{g}}(\mathfrak{aff}^*_\mathfrak{g})$.

\subsection*{Contact structures} The following example involves a groupoid which is non-flat. Recall that a (co-orientable) contact structure on a $2k+1$-manifold $B$ -- i.e. $\alpha\in \Omega^1(B)$ such that $\alpha\wedge (d\alpha)^n$ -- is in particular a $\mathrm{Sp}(k,1)$-structure, where $\mathrm{Sp}(k,1)\subset \mathrm{Gl}_{2k+1}(\mathbb{R})$ consist of block matrices
\[
\begin{pmatrix}
A & v \\
0 & \lambda 
\end{pmatrix},\quad A\in \mathrm{Sp}(k),\ \lambda > 0, v\in \mathbb{R}^{2k}.
\]
Let us denote $(P^\mathrm{\alpha^\mathrm{std}},\tau)$ the $\mathrm{Sp}(k,1)$-structure corresponding to the standard contact form on $\mathbb{R}^{2k+1}$ with coordinates $(x^1,\dots, x^k,y,z^1,\dots z^k)$, i.e. the form
\[
\alpha^\mathrm{std} = dy-\sum\limits_{i=1}^k z^idx^i.
\]
It is not difficult to see that $(P^\mathrm{\alpha^\mathrm{std}},\tau)$ is not isomorphic to $(P^\textrm{$G$-flat},\tau)$; indeed, the latter corresponds to the codimension one foliation $\F^\mathrm{flat}$ tangent to the kernel of 
\[
dy = 0
\]
and equipped with the leafwise symplectic form
\[
\omega^{\F^\mathrm{flat}} = \sum\limits_{i=1}^k dx^i\wedge dz^i.
\]
In fact, $(P^\mathrm{\alpha^\mathrm{std}},\tau)$ is isomorphic to 
\[
(\mathbb{R}^{2k+1} \times \mathrm{Sp}(k,1), A\cdot \pr_{\mathbb{R}^{2k+1}}^*\omega^\mathrm{MC}_{\mathrm{Hei}})
\]
where $\mathrm{Hei}$ is $\mathbb{R}^{2k+1}$ with the Heisenberg group structure and the notation $A\cdot \pr_{\mathbb{R}^{2k+1}}^*\omega^\mathrm{MC}_{\mathrm{Hei}}$ is the same we used in the previous example. The isomorphism leads us to the $\mathrm{Sp}(k,1)$-structure with connection
\[
(\mathbb{R}^{2k+1} \times \mathrm{Sp}(k,1), \theta^{\alpha^\mathrm{std}})
\]
where this time
\[
\theta^{\alpha^\mathrm{std}} = \pr_{\mathbb{R}^{2k+1}}^*\omega^\mathrm{MC}_{\mathrm{Hei}} + \pr_{\mathrm{Sp}(k,1)}^*\omega^\mathrm{MC}_{\mathrm{Sp}(k,1)}
\]
Differently form the case of $(P^\textrm{$G$-flat},\theta^\textrm{$G$-flat})$, $(P^\mathrm{\alpha^\mathrm{std}},\theta^{\alpha^\mathrm{std}})$ is not a flat Cartan geometry. The resulting gauge groupoid
\[
(\G^{\alpha^\mathrm{std}},\Theta^{\alpha^\mathrm{std}})
\]
takes the form
\[
\G^{\alpha^\mathrm{std}} = \mathbb{R}^n\times \mathrm{Sp}(k,1)\times \mathbb{R}^n
\]
\[
\Theta^{\alpha^\mathrm{std}} = \pr_{\mathrm{Sp}(k,1)}^*\omega^\mathrm{MC}_{\mathrm{Sp}(k,1)}+\delta\omega^\mathrm{MC}_{\mathrm{Hei}}.
\]
This is not a flat groupoid. 

Any $\mathrm{Sp}(k,1)$-structure with connection $(P,\theta)$ can be presented either as a principal $(\G^{\textrm{$G$-flat}},\Theta^{\textrm{$G$-flat}})$-bundle or as a principal $(\G^{\alpha^\mathrm{std}},\Theta^{\alpha^\mathrm{std}})$. The two presentations produce different invariants. While the characteristic classes that arise from the $(\G^{\textrm{$G$-flat}},\Theta^{\textrm{$G$-flat}})$ action measure obstruction to flatness -- recall Theorem \ref{thm:chern_weil_CG} -- the ones arising from $(\G^{\alpha^\mathrm{std}},\Theta^{\alpha^\mathrm{std}})$ via our Theorem \ref{thm:main_thm} should be interpreted as obstructions for $(P,\theta)$ being isomorphic to $(P^\mathrm{\alpha^\mathrm{std}},\theta^{\alpha^\mathrm{std}})$.

\end{document}